\documentclass[11pt]{article}
\usepackage{amsthm,amssymb}
\usepackage{graphicx,subfigure}
\usepackage{amsmath}    % need for subequations
\usepackage{graphicx}   % need for figures
\usepackage{verbatim}   % useful for program listings
\usepackage{color}      % use if color is used in text
\usepackage{subfigure}  % use for side-by-side figures
\usepackage{hyperref}   % use for hypertext links, including those to external documents and URLs
\newcommand{\del}{\backslash}

 \usepackage[usenames,dvipsnames]{pstricks}
  \usepackage{epsfig,pst-node,pst-tree,pst-grad}
\usepackage{lineno}
\newtheorem{theorem}{Theorem}[section]

\newtheorem{lemma}[theorem]{Lemma}

\newtheorem{corollary}[theorem]{Corollary}

\newcommand{\be}{\begin{enumerate}}

\newcommand{\ee}{\end{enumerate}}

\newcommand{\su}{\subseteq}

\newcommand{\noi}{\noindent}

\begin{document}

\title{Characterizing binary matroids with no  $P_9$-minor}

\author{{Guoli Ding$^{1}$ ~and Haidong Wu$^{2}$}\\
{\small 1. Department of Mathematics,  Louisiana State University,
Baton Rouge, Louisiana, USA}\\
{\small Email: ding@math.lsu.edu}\\
{\small 2.  Department of Mathematics,  University of Mississippi, University,
Mississippi, USA} \\
{\small Email: hwu@olemiss.edu}}

\date{}
\maketitle

\begin{abstract}
 In this paper, we give a complete characterization of binary matroids with no $P_9$-minor. A 3-connected binary matroid $M$ has no $P_9$-minor  if and only if $M$ is
 one of the internally 4-connected  non-regular 
 minors of a special 16-element matroid $Y_{16}$,   a 3-connected regular matroid, a binary spike with rank at least four,  or  a matroid obtained by 
 3-summing copies of the Fano matroid to a 3-connected cographic matroid $M^*(K_{3, n})$, $M^*(K_{3, n}^{\prime})$, $M^*(K_{3, n}^{\prime\prime})$, or 
$M^*(K_{3, n}^{\prime\prime\prime})$ ($n\ge 2$). Here the simple graphs $K_{3, n}^{\prime}, K_{3, n}^{\prime\prime}$, and 
$K_{3, n}^{\prime\prime\prime}$ are obtained  from $K_{3, n}$ by adding one, two, or  three edges in the color class of size three, respectively.

\end{abstract}

\section{Introduction}\label{intro}
It is well known that the class of binary matroids consists of all matroids without any $U_{2, 4}$-minor, and the class of regular matroids consists of matroids without any $U_{2, 4}, F_7$ or $F_7^*$-minor.   Kuratowski's Theorem states that a graph is planar if and only if it has no minor that is isomorphic to    $K_{3,3}$ or $K_5$.  These examples show that characterizing a class of graphs and matroids without certain minors is often of fundamental importance. We say that a matroid is {\it $N$-free} if it does not contain a minor that is isomorphic to $N$.     A $3$-connected matroid $M$ is said to be internally $4$-connected if  for any $3$-separation of $M$, one side of the separation is either a triangle or a triad.

There is much interest in characterizing binary matroids without small 3-connected minors. Since non-3-connected matroids can be constructed by 3-connected matroids using 1-, 2-sum operations, one needs only determine the 3-connected members of a minor closed class. There is exactly one 3-connected binary matroid with 6-elements, namely, $W_3$ where $W_n$ denotes both the   
 wheel graph with $n$-spokes and the cycle matroid of $W_n$. 
 There are exactly two 7-element binary 3-connected matroids, $F_7$ and $F_7^*$.  There are three 8-element binary 3-connected matroids, $W_4$, $S_8$ and $AG(3,2)$, and there are eight  9-element 3-connected binary matroids: $M(K_{3,3})$, $M^*(K_{3,3})$, Prism, $M(K_5\backslash e)$, $P_9, P_9^*$,  binary spike $Z_4$ and its dual $Z_4^*$. 

  \vspace{0.2in}
  \begin{tabular}{|r|r|r|}
\hline
$|E(M)|$& Binary 3-connected matroids\\
\hline
6 & $W_3$\\
7 & $F_7, F_7^*$ \\
8& $W_4, S_8, AG(3,2)$\\
9& $M(K_{3,3}),  M^*(K_{3,3}), M(K_5\backslash e), Prism, P_9, P_9^*, Z_4, Z_4^*$\\
\hline
\end{tabular}

 \vspace{0.2in}
   
\begin{figure}   
\begin{center}
\includegraphics[scale=0.6]{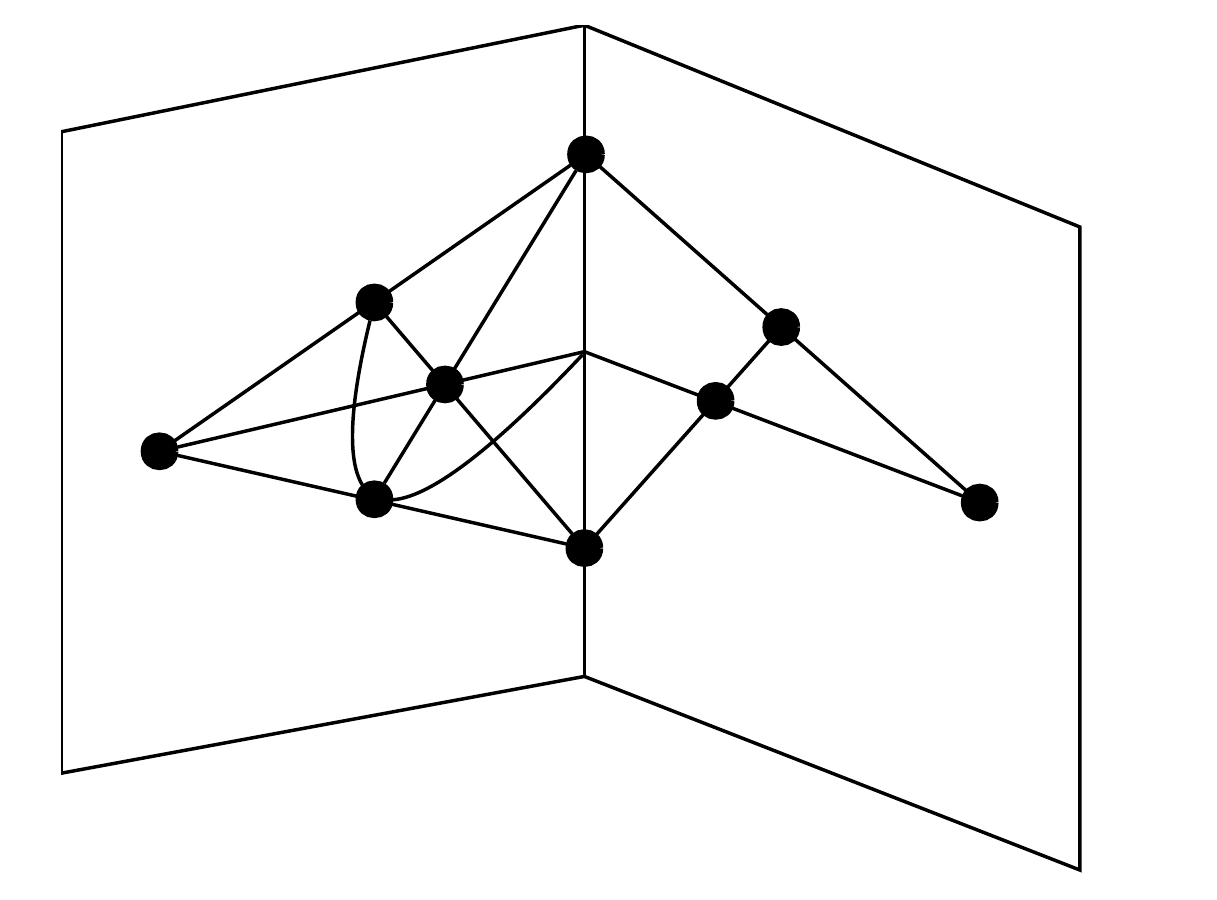}
\caption{\hspace{0.1in} A geometric representation of $P_9$}
%\label{$P_9$}
\end{center}
\end{figure}

For each matroid $N$ in the above list with less than nine elements, with the exception of $AG(3, 2)$, 
 the problem of characterizing 3-connected binary matroids with no  $N$-minor has been solved. Since every 3-connected binary matroid having at least four elements has a 
 $W_3$-minor, the class of 3-connected binary matroids excluding $W_3$ contains only the trivial 3-connected matroids with at most three elements. Seymour in \cite{seymour1980} determined all 3-connected binary matroids with no $F_7$-minor ($F_7^*$-minor).  Any such matroid is either regular or is isomorphic to $F_7^*$ ($F_7$).  
 In \cite{oxley87}, Oxley characterized all 3-connected binary  $W_4$-free matroids. These are exactly $M(K_4)$, $F_7, F_7^*$, binary spikes $Z_r$, $Z_r^*$, $Z_r\backslash t$, or $Z_r\backslash y_r$ ($r\ge 4$) plus the trivial 3-connected matroids with at most three elements. 
 It is well known that $F_7, F_7^*$, and $AG(3, 2)$ are the only 3-connected binary non-regular matroids without any $S_8$-minor. 
 
 In the book \cite{mrw}, Mayhew,  Royle and Whittle characterized all internally 4-connected binary $M(K_{3,3})$-free matroids.   Mayhew and Royle \cite{mayhew2012}, and independently Kingan and Lemos \cite{kingan-lemos}, determined all  internally $4$-connected binary Prism-free (therefore $M(K_5\backslash e)$-free) matroids.  For each matroid $N$ in the above list with exactly  nine elements, 
 the problem of characterizing 3-connected binary matroids with no $N$-minor is still unsolved yet.  The problem of characterizing internally 4-connected binary $AG(3,2)$-free matroids is also open.    Since $Z_4$ has an $AG(3,2)$-minor, characterizing internally 4-connected binary $Z_4$-free matroids is an even harder problem.   Oxley \cite{oxley87} determined all 3-connected binary matroids with no $P_9$- or $P_9^*$-minor:

\begin{theorem} \label{oxleyp9}
Let $M$ be a binary matroid. Then $M$ is 3-connected having no minor isomorphic to $P_9$ or $P_9^*$ if and only if 

(i) $M$ is regular and 3-connected;

(ii) $M$ is a binary spike $ Z_r, Z_r^*, Z_r\backslash y_r$ or $Z_r\backslash t$ for some $r\ge 4$; or

(iii) $M\cong F_7$ or $F_7^*$.

\end{theorem}

$P_9$ is  a very important matroid and it appears frequently in the structural matroid theory (see, for example, \cite{moow2011, oxley87, zhou2004}). In this paper, we give a complete characterization of  the 3-connected binary matroids with no  $P_9$-minor. Before we state our main result,  we describe a class of 
non-regular matroids. First let $\cal K$ be the class 3-connected cographic matroids $N=M^*(K_{3, n})$, $M^*(K_{3, n}^{\prime})$, $M^*(K_{3, n}^{\prime\prime})$, or 
$M^*(K_{3, n}^{\prime\prime\prime})$ ($n\ge 2$). Here the simple graphs $K_{3, n}^{\prime}, K_{3, n}^{\prime\prime}$, and 
$K_{3, n}^{\prime\prime\prime}$ are obtained 
from $K_{3, n}$ by adding one, two, or  three edges in the color class of size three, respectively. 
Note that when $n=2$, $N\cong W_4, $ or the cycle matroid of 
the prism graph. From now on, we will use Prism to denote the prism graph as well as its cycle matroid. Take any 
$t$ disjoint triangles $T_1, T_2, \dots, T_t$ ($1\le t\le n$) of $N$ and $t$ copies of $F_7$. Perform 3-sum operations consecutively starting from $N$ and $F_7$ along the triangles $T_i$ ($1\le i\le t$).  Any resulting matroid in this infinite class of matroids is called a (multi-legged) starfish. Note that each starfish is not regular since at least one Fano was used (and therefore the resulting matroid has an $F_7$-minor) in the construction. The class of starfishes and the class of spikes have empty intersection as spikes are $W_4$-free, while each
starfish has a $W_4$-minor.

Our next result, the main result of this paper, generalizes Oxley's Theorem \ref{oxleyp9} and completely determines the 3-connected $P_9$-free binary matroids. The matroid $Y_{16}$, a single-element extension of $PG(3,2)^*$,  in standard representation without the identity matrix is given in Figure 2.

\begin{theorem}    \label{main}
Let M be a binary matroid. Then $M$ is 3-connected  having no minor isomorphic to $P_9$ if and only if one of the following is true:
\be
\item[(i)] $M$ is one of the $16$  internally 4-connected non-regular minors of $Y_{16}$; or
\item[(ii)] $M$ is regular and 3-connected; or 
\item[(iii)] $M$ is a binary spike $ Z_r, Z_r^*, Z_r\backslash y_r$ or $Z_r\backslash t$ for some $r\ge 4$; or 
\item [(iv)] $M$  is  a   starfish. 
\ee
\end{theorem}   

\begin{figure}[!h]
\begin{center}
$ \left[ \begin{array}{ccccc}
1 & 1 & 1 & 0 & 0  \\
1 & 1 & 0 & 1 & 0  \\
1 & 0 & 1 & 1 & 0 \\
0 & 0 & 1 & 1 & 1  \\
0 & 1 & 0 & 1 & 1 \\
0 & 1 & 1 & 0 & 1 \\
0 & 1 & 1 & 1 & 0 \\
1 & 0 & 0 & 1 & 1 \\
1 & 0 & 1 & 0 & 1 \\
1 & 1 & 0 & 0 & 1 \\
1 & 1 & 1 & 1 & 1 \\
 
\end{array} \right]$
\caption{\hspace{0.1in}A binary standard representation for $Y_{16}$ }
\label{matroidp17}
\end{center}
\end{figure}

The next result, which follows easily from the last theorem, characterizes all binary $P_9$-free matroids. 

\begin{theorem}    \label{main2}
Let M be a binary matroid. Then $M$ has no minor isomorphic to $P_9$ if and only if $M$ can be constructed from internally 4-connected non-regular minors of 
$Y_{16}$, 3-connected regular matroids, binary spikes, and starfishes using the 
operations of direct sum and 2-sum. 
\end{theorem}   

\proof  Since every matroid can be constructed from 3-connected proper minors of itself by the operations of direct sum and 2-sum, by Theorem \ref{main}, 
the forward direction is true. Conversely, suppose that $M=M_1\oplus M_2$, or $M=M_1\oplus_2 M_2$, where $M_1$ and $M_2$ are both $P_9$-free. As $P_9$ is 3-connected, by \cite[Proposition 8.3.5]{oxley}, $M$ is also $P_9$-free. Thus 
if $M$ is  constructed from internally 4-connected non-regular minors of 
$Y_{16}$, 3-connected regular matroids, binary spikes, and starfishes using the 
operations of direct sum and 2-sum, then $M$ is also $P_9$-free. \qed

Our proof does not use Theorem 1.1 except we use the fact that all spikes are $P_9$-free which can be proved by an easy induction argument.  In Section 2, we determine all internally 4-connected binary $P_9$-free matroids.  These are exactly the $16$  internally 4-connected non-regular minors of $Y_{16}$.  These matroids are 
determined using the Sage matroid package and the computation is confirmed by the matroid software Macek.     Most of the work is in  Section 3, which  is to determine how the internally 4-connected pieces can be put together to avoid a $P_9$-minor.  

For terminology we follow \cite{oxley}.  Let $M$ be a matroid.  The {\it connectivity function} $\lambda_M$ of $M$  is defined as follows. For $X \su E$ let
\begin{equation}
\lambda_M(X)=r_M(X)+r_M(E-X)-r(M).
\end{equation}
Let $k \in \mathbb{Z}^+$. Then both $X$ and $E-X$ are said to be {\it k-separating}  if $\lambda_M(X)=\lambda_M(E-X) < k$. If $X$ and $E-X$ are $k$-separating and $\min \{|X|,|E-X|\} \geq k$, then $(X,E-X)$ is said to be a $k$-separation of $M$.
Let $\tau (M)=\min \{~j :M ~\text{has a } j\text{-separation}\}$ if $M$ has a $k$-separation for some $k$; otherwise let $\tau (M) =\infty.$   $M$ is $k$-{\it connected} if $\tau (M) \geq k$.  
Let $(X,E-X)$ be a $k$-separation of $M$. This separation is said to be a {\it minimal k-separation} if
$\min \{|X|,~|E-X|\}=k$. The matroid $M$ is called {\it internally 4-connected} if and only if $M$ is $3$-connected and the only $3$-separations of $M$ are minimal (in other words, either $X$ or $Y$ is a triangle or a triad). %$M$ is {\it weakly 4-connected} if $M$ is 3-connected, and for any 3-separation $(X, Y)$ of $M$, either $|X|\le 4$ or $|Y|\le 4$. 

\section{Characterizing  internally 4-connected binary $P_9$-free matroids}

In this section, we determine all  internally 4-connected  binary $P_9$-free matroids.

\begin{theorem} \label{maini4c}
A binary matroid $M$ is internally 4-connected  and $P_9$-free if and only if 

(i) $M$ is internally 4-connected graphic or cographic; or

(ii) $M$ is one of the 16 internally 4-connected non-regular minors of $Y_{16}$; or

(iii) $M$ is isomorphic to $R_{10}$. 

\end{theorem}

Sandra Kingan recently informed us that she also obtained the internally 4-connected binary $P_9$-free matroids as a consequence of 
a decomposition result for 3-connected binary $P_9$-free matroids.

The following two well-known theorems of Seymour  \cite{seymour1980} will be used  in our proof. 

\begin{theorem} {(Seymour's Splitter Theorem)} \label{SST}
 Let $N$ be a $3$-connected proper minor of a 3-connected matroid $M$  such that $|E(N)|\ge 4$ and if $N$ is a wheel, it is the largest wheel minor of $M$; while if $N$ is a whirl,  it is the largest whirl minor of $M$. Then $M$ has a 3-connected minor $M^{\prime}$ which is isomorphic to a single-element extension or coextension of  $N$.
\end{theorem}

\begin{theorem}\label{seymour-3-conn}
If $M$ is an internally 4-connected regular matroid, then M is graphic, cographic, or is isomorphic to $R_{10}$. 
\end{theorem}

 The following result is due to  Zhou \cite[Corollary 1.2]{zhou2004}.
 
\begin{theorem}\label{zhou2004}

A non-regular internally 4-connected binary matroid other than $F_7$ and $F_7^*$ contains one of the following matroids as a minor: $N_{10}$, 
$\widetilde {K_5}$, $\widetilde {K_5}^*$, $T_{12}\backslash e$, and $T_{12}/e$. 

\end{theorem}

The matrix representations of these matroids  can be found in \cite{zhou2004}. We use $X_{10}$ to denote  the matroid $\widetilde {K_5}^*$. It is 
straightforward to verify that among the five matroids in Theorem \ref{zhou2004}, only $X_{10}$ has no $P_9$-minor. 
We use $\cal L$ to denote the set of matroids consisting of the  following matroids in reduced standard representation,  in addition to  
$F_7, F_7^*$ and $Y_{16}$. 
From the matrix representations of these matroids,  it is straightforward to check that each matroid in $\cal L$ is a minor of $Y_{16}$, and 
each has an $X_{10}$-minor. Indeed, It is clear that  (i) each $X_i$ is a single-element co-extension of $X_{i-1}$ 
 for $11\le i\le 15$; (ii) each $Y_i$ is a single-element extension of $X_{i-1}$ for $11\le i\le 16$; (iii) each $Y_i$ is a single-element co-extension of 
 $Y_{i-1}$ for $11\le i\le 16$, and it is easy to check that (iv)  in the list $X_{10}, X'_{11},  X'_{12}, X_{13}$, each matroid is a single-element coextension of its
  immediate predecessor. 
Therefore,  $X_{10}$ is a minor of all matroids in $\cal L$, and each is a minor of $Y_{16}$.  From these matrices, 
 it is also routine to check that the only matroid of $ {\cal L}$ having a triangle is $F_7$ (this can also be easily 
verified by using the Sage matroid package).

\[ X_{10}: \left( \begin{array}{cccc}
1 & 1 & 1 & 0   \\
1 & 1 & 0 & 1   \\
1 & 0 & 1 & 1  \\
0 & 0 & 1 & 1   \\
0 & 1 & 0 & 1  \\
1 & 0 & 0 & 1 

\end{array} \right) X_{11}: \left( \begin{array}{cccc}
1 & 1 & 1 & 0   \\
1 & 1 & 0 & 1   \\
1 & 0 & 1 & 1  \\
0 & 0 & 1 & 1   \\
0 & 1 & 0 & 1  \\
0 & 1 & 1 & 0\\
1 & 0 & 0 & 1 

\end{array} \right) X_{11}^{\prime}: \left( \begin{array}{ccccc}
1 & 1 & 1 & 0   \\
1 & 1 & 0 & 1   \\
1 & 0 & 1 & 1  \\
0 & 0 & 1 & 1   \\
0 & 1 & 0 & 1  \\
0 & 1 & 1 & 1\\
1 & 0 & 0 & 1 

\end{array} \right) Y_{11}: \left( \begin{array}{ccccc}
1 & 1 & 1 & 0  & 0 \\
1 & 1 & 0 & 1  & 0 \\
1 & 0 & 1 & 1  &0\\
0 & 0 & 1 & 1   &1\\
0 & 1 & 0 & 1  & 1\\
1 & 0 & 0 & 1  & 1

\end{array} \right)\]

\[ X_{12}: \left( \begin{array}{cccc}
1 & 1 & 1 & 0   \\
1 & 1 & 0 & 1   \\
1 & 0 & 1 & 1  \\
0 & 0 & 1 & 1   \\
0 & 1 & 0 & 1  \\
0 & 1 & 1 & 0 \\
0 & 1 & 1 & 1 \\
1 & 0 & 0 & 1 

\end{array} \right) X_{12}^{\prime}: \left( \begin{array}{ccccc}
1 & 1 & 1 & 0   \\
1 & 1 & 0 & 1   \\
1 & 0 & 1 & 1  \\
0 & 0 & 1 & 1   \\
0 & 1 & 0 & 1  \\
0 & 1 & 1 & 1  \\ 
1 & 0 & 0 & 1 \\
1 & 1& 1 & 1

\end{array} \right) Y_{12}: \left( \begin{array}{ccccc}
1 & 1 & 1 & 0  & 0 \\
1 & 1 & 0 & 1  & 0 \\
1 & 0 & 1 & 1  &0\\
0 & 0 & 1 & 1   &1\\
0 & 1 & 0 & 1  & 1\\
0 & 1 & 1 & 0  & 1\\
1 & 0 & 0 & 1  & 1

\end{array} \right)\]

\[ X_{13}: \left( \begin{array}{cccc}
1 & 1 & 1 & 0   \\
1 & 1 & 0 & 1   \\
1 & 0 & 1 & 1  \\
0 & 0 & 1 & 1   \\
0 & 1 & 0 & 1  \\
0 & 1 & 1 & 0  \\
0 & 1 & 1 & 1   \\
1 & 0 & 0 & 1 \\
1 & 0 & 1& 0
\end{array} \right) Y_{13}: \left( \begin{array}{ccccc}
1 & 1 & 1 & 0  & 0 \\
1 & 1 & 0 & 1  & 0 \\
1 & 0 & 1 & 1  &0\\
0 & 0 & 1 & 1   &1\\
0 & 1 & 0 & 1  & 1\\
0 & 1 & 1 & 0   &1\\
0 & 1 & 1 & 1  & 0\\
1 & 0 & 0 & 1  & 1
\end{array} \right)  X_{14}: \left( \begin{array}{cccc}
1 & 1 & 1 & 0   \\
1 & 1 & 0 & 1   \\
1 & 0 & 1 & 1  \\
0 & 0 & 1 & 1   \\
0 & 1 & 0 & 1  \\
0 & 1 & 1 & 0  \\
0 & 1 & 1 & 1   \\
1 & 0 & 0 & 1 \\
1 & 0 & 1 & 0\\
1 & 1 & 0& 0
\end{array} \right)  Y_{14}: \left( \begin{array}{ccccc}
1 & 1 & 1 & 0  & 0 \\
1 & 1 & 0 & 1  & 0 \\
1 & 0 & 1 & 1  &0\\
0 & 0 & 1 & 1   &1\\
0 & 1 & 0 & 1  & 1\\
0 & 1 & 1 & 0   &1\\
0 & 1 & 1 & 1  & 0\\
1 & 0 & 0 & 1  & 1\\
1 & 0 & 1 & 0 & 1

\end{array} \right)\]

\begin{align}
X_{15}\cong PG(3,2)^* &:\begin{pmatrix}
1 & 1 & 1 & 0   \\
1 & 1 & 0 & 1   \\
1 & 0 & 1 & 1  \\
0 & 0 & 1 & 1   \\
0 & 1 & 0 & 1  \\
0 & 1 & 1 & 0 \\
0 & 1 & 1 & 1  \\
1 & 0 & 0 & 1 \\
1 & 0 & 1 & 0 \\
1 & 1 & 0 & 0 \\
1 & 1 & 1 & 1 
\end{pmatrix}
&
Y_{15} &:\begin{pmatrix}
1 & 1 & 1 & 0 & 0  \\
1 & 1 & 0 & 1 & 0  \\
1 & 0 & 1 & 1 & 0 \\
0 & 0 & 1 & 1 & 1  \\
0 & 1 & 0 & 1 & 1 \\
0 & 1 & 1 & 0 & 1 \\
0 & 1 & 1 & 1 & 0 \\
1 & 0 & 0 & 1 & 1 \\
1 & 0 & 1 & 0 & 1 \\
1 & 1 & 0 & 0 & 1  
\end{pmatrix}
\end{align}

%It is easy to see from the matrix representations that  (i) each $X_i$ is a single-element co-extension of $X_{i-1}$ for $11\le i\le 15$; (ii) each $Y_i$ is a single-%element extension of $X_{i-1}$ for $11\le i\le 16$; (iii) each $Y_i$ is a single-element co-extension of $Y_{i-1}$ for $11\le i\le 16$, and (iv) $X_{10}$ is a minor of all %matroids listed above, and each is a minor of $Y_{16}$. Note that $X^*_{15}\cong PG(3,2)$, so $Y_{16}$ is a single-element extension of  $PG(3,2)^*$.

\noindent {\it Proof of Theorem \ref{maini4c}:  } If $M$ is one of the matroids listed in (i) to (iii), then  $M$ is internally 4-connected.  All matroids in (i) or (iii) are regular, thus are $P_9$-free. 
Using the Sage matroid package,  it is easy to verify that $Y_{16}$ is $P_9$-free, hence all matroids in (ii) are also $P_9$-free. Let $M$ be an internally  4-connected binary matroid with no $P_9$-minor. If $M$ is regular, then by Theorem \ref{seymour-3-conn}, $M$ is either graphic, cographic, or isomorphic to $R_{10}$, which is regular. Therefore, we need only show that  an internally 4-connected matroid $M$ is  non-regular and $P_9$-free if and only if $M$ is a non-regular minor of $Y_{16}$. Suppose that $M$ is  an internally 4-connected 
non-regular and $P_9$-free matroid. If $M$ has exactly seven elements, then
 $M\cong F_7$ or $M\cong F_7^*$.  Suppose that $M$ has at least eight elements. By Theorem \ref{zhou2004},   $M$ has an $N_{10}$, $X_{10}$, $X_{10}^*$, $T_{12}\backslash e$, or $T_{12}/e$-minor. Since all but $X_{10}$ has a $P_9$-minor among these matroids, $M$ must have an $X_{10}$-minor. We use  the Sage matroid package (by writing simple Python scripts) and the matroid software Macek  independently to  do our computation and have obtained the same result. 
Excluding $P_9$, we extend  and coextend $X_{10}$ seven times and found only  thirteen 3-connected binary matroids. These matroids are $X_{11}$, $X_{11}^{\prime}, Y_{11}, X_{12}, X_{12}^{\prime}, Y_{12}, X_{13}, Y_{13},  X_{14}, Y_{14}, X_{15}\cong PG(3,2)^*, Y_{15},$ and $Y_{16}$; each having at most 16 elements; each being a minor of $Y_{16}$;  and 
each being internally 4-connected. As $X_{10}$ is neither a wheel nor a whirl, by the Splitter Theorem (Theorem \ref{SST}), $M$ is one of the matroids in $\cal L$, each of which is a non-regular internally 4-connected minor of $Y_{16}$. Note that all non-regular internally 4-connected minors of $Y_{16}$  are $P_9$-free,  hence $\cal L$ consists of all internally 4-connected non-regular minors of $Y_{16}$. \qed

\section{Characterizing 3-connected binary $P_9$-free matroids}

In this section, we will prove our main result. We begin with several lemmas. Let $G$ be a graph with a specified triangle $T=\{e_1,e_2,e_3\}$. By a {\it
rooted $K_4''$-minor} using $T$ we mean a loopless minor $H$ of $G$ such that $si(H)\cong K_4$;
$\{e_1,e_2,e_3\}$ remains a triangle of $H$; and $H\backslash \{e_i,e_j\}$
is isomorphic to $K_4$, for some distinct $i,j\in\{1,2,3\}$. By a {\it
rooted $K_4'$-minor} using $T$ we mean a loopless minor $H$ of $G$ such that $si(H)\cong K_4$;
$\{e_1,e_2,e_3\}$ remains a triangle of $H$; and $H\backslash e_i$
is isomorphic to $K_4$, for some $i\in\{1,2,3\}$. 
Let $T$ be a specified triangle of a matroid $M$. We can define a rooted $M(K_4')$-minor using $T$ and 
a rooted $M(K_4'')$-minor using $T$ similarly. Moreover, in the following proof, any $K_4^\prime$ is obtained from 
$K_4$ by adding a parallel edge to an element in the common triangle $T$ used in the 3-sum specified in the context.

\begin{lemma} (\cite{seymour1985}) \label{k4} Let $T$ be a triangle of 3-connected binary matroid $M$ with at least four elements. Then $T$ is contained in a $M(K_4)$-minor of $M$. 
\end{lemma}

\begin{lemma} (\cite{seymour1984}) \label{seymour-F7} Let $T$ be a triangle of a binary non-graphic matroid $M$. Then the following are true:

(i) If $M$ is non-regular, then $T$ is contained in a $F_7$-minor; 

(ii) If $M$ is regular but not graphic, then $T$ is contained in a $M^*(K_{3,3})$-minor.

\end{lemma}

Let $M_1$ and $M_2$ be matroids with ground sets $E_1$ and $E_2$ such that $E_1\cap E_2=T$ and $M_1|T=M_2|T=N$. The following result of Brylawski \cite{brylawski} about the generalized parallel connection can be found in \cite[Propsition 11.4.14]{oxley}.

\begin{lemma} \label{brylawski}
The generalized parallel connection $P_N(M_1,M_2)$ has the following properties:
\be
\item[(i)] $P_N(M_1,M_2)|E_1=M_1$ and $P_N(M_1,M_2)|E_2=M_2$.
\item[(ii)] If $e\in E_1 - T$, then $P_N(M_1, M_2)\del e=P_N(M_1\del e, M_2)$.
\item[(iii)] If $e\in E_1 - cl_1(T)$, then $P_N(M_1 ,M_2)/e=P_N(M_1/e, M_2)$.
\item[(iv)] If $e\in E_2 - T$, then $P_N(M_1, M_2)\del e=P_N(M_1, M_2\del e)$.
\item[(v)] If $e\in E_2 - cl_2(T)$, then $P_N(M_1 ,M_2)/e=P_N(M_1, M_2/e)$.
\item[(vi)] If $e\in T$, then $P_N(M_1, M_2)/e=P_{N/e}(M_1/e, M_2/e)$.
\item[(vii)] $P_N(M_1, M_2)/T=(M_1/T)\oplus (M_2/T)$.
\ee
\end{lemma}

In the rest of this paper, we consider the case when the generalized parallel connection is defined across a triangle $T$, where $T$ is the common triangle of 
the binary matroids  $M_1$ and $M_2$. Then $P_N(M_1, M_2)=P_N(M_2, M_1)$ 
(see  \cite[Propsition 11.4.14]{oxley}). Moreover, $N=M_1|T=M_2|T\cong U_{2,3}$.  We will use $T$ to denote both the triangle and the submatroid $M_1|T$. 
Thus we use $P_T(M_1, M_2)$ instead of $P_N(M_1, M_2)$ for the rest of the paper. 

\begin{lemma}\label{order}

Let $M=P_T(M_1, P_S(M_2, M_3))$ where $M_i$ is a  binary matroid ($1\le i\le 3$); $S$ is the common triangle of $M_2$ and $M_3$; $T$ is the common triangle of $M_1$ and $M_2$.  Then the following are true:

(i) if  $E(M_1)\cap (E(M_3)\backslash E(M_2))=\emptyset$, then  $M= P_S(P_T(M_1, M_2), M_3)$; 

(ii) if $E(M_1)\cap E(M_3)=\emptyset$, then $M_1\oplus_3 (M_2\oplus_3 M_3)=(M_1\oplus_3 M_2)\oplus_3 M_3$.

\end{lemma}

\proof (i) As $E(M_1)\cap (E(M_3)\backslash E(M_2))=\emptyset$, $T=E(M_1)\cap E(P_S(M_2, M_3))$, and $T$ is the common triangle of $M_1$ and $P_S(M_2, M_3)$. 
 Moreover, $S=E(M_3)\cap E(P_T(M_1, M_2))$, and $S$ is the common triangle of $M_3$ and $P_T(M_1, M_2)$. By  \cite[ Proposition 11.4.13]{oxley}, a set $F$ of $M$ is a flat if and only if $F\cap E(M_1)$ is a flat of $M_1$ and $F\cap E(P_S(M_2, M_3))$ is a flat of 
$P_S(M_2, M_3)$. The latter is true if and only if $[F\cap (E(M_2)\cup E(M_3))]\cap E(M_i)=F\cap E(M_i)$ is a flat of $M_i$ for $i=2, 3$. Therefore, $F$ is a 
flat of $M$ if and only if $F\cap E(M_i)$ is a flat of $M_i$ for $1\le i\le 3$. The same holds for $P_S(P_T(M_1, M_2), M_3)$. Thus $M=P_S(P_T(M_1, M_2), M_3)$. 

(ii) As $E(M_1)\cap E(M_3)=\emptyset$,  we deduce that $S\cap T=\emptyset$, and the conclusion of (i) holds. Therefore,
$$P_T(M_1, P_S(M_2, M_3))\backslash (S\cup T)= P_S(P_T(M_1, M_2), M_3)\backslash (S\cup T).$$  By Lemma \ref{brylawski}, we conclude that 
$$P_T(M_1, P_S(M_2, M_3)\backslash S)\backslash  T= P_S(P_T(M_1, M_2)\backslash T, M_3)\backslash S.$$ That is, $M_1\oplus_3 (M_2\oplus_3 M_3)=(M_1\oplus_3 M_2)\oplus_3 M_3$. \qed

\begin{lemma}\label{cocircuit}

Let $M=P_T(M_1, M_2)$ where $M_i$ is a  binary matroid ($1\le i\le 2$) and $T$ is the common triangle of $M_1$ and $M_2$. Then $C^*$ is a cocircuit of $M$ if and only if one of the following is true:

(i) $C^*$ is a cocircuit of $M_1$ or $M_2$ avoiding $T$;

(ii) $C^*=C_1^*\cup C_2^*$ where $C_i^*$ is a cocircuit of $M_i$ such that $C_1^*\cap T=C_2^*\cap T$, which has exactly two elements.  
\end{lemma}

\proof By  \cite[ Proposition 11.4.13]{oxley}, a set $F$ of $M$ is a flat if and only if $F\cap E(M_i)$ is a flat of $M_i$ for $1\le i\le 2$. Moreover, for any flat $F$ of $M$, $r(F)=r(F\cap E(M_1))+r(F\cap E(M_2))-r(F\cap T)$ (see, for example,  \cite[(11.23)]{oxley}). Let $C^*$ be a cocircuit of $M$ and $H=E(M)-C^*$.  As $M$ is binary, $|C^*\cap T|=0, 2$, and thus  $|H\cap T|=3, 1$. First assume that $|C^*\cap T|=0$. As $r(H)=r(H\cap E(M_1))+r(H\cap E(M_2))-r(H\cap T)$, then $r(M)-1=r(M_1)+r(M_2)-3=r(H)=r(H\cap E(M_1))+r(H\cap E(M_2))-2$. Thus,
$$ r(M_1)+r(M_2)-1=r(H\cap E(M_1))+r(H\cap E(M_2)).$$

Therefore, either $r(H\cap E(M_1))=r(M_1)-1$ and $r(H\cap E(M_2))=r(M_2)$, or $r(H\cap E(M_2))=r(M_2)-1$ and $r(H\cap E(M_1))=r(M_1)$.  In the former case, as  $H\cap E(M_1)$ and $H\cap E(M_2)$ are flats of $M_1$ and $M_2$ respectively, we deduce that $H\cap E(M_2)=E(M_2)$;  $H\cap E(M_1)$ is a hyperplane of $M_1$ and thus $C^*\subseteq E(M_1)$ is a cocircuit of $M_1$ avoiding $T$. The latter case is similar. 

If $|C^*\cap T|=2$, then $|H\cap T|=1$. As $r(H)=r(H\cap E(M_1))+r(H\cap E(M_2))-r(H\cap T)$, we deduce that  $r(M)-1=r(M_1)+r(M_2)-3=r(H)=r(H\cap E(M_1))+r(H\cap E(M_2))-1$. We conclude that
$$ r(M_1)+r(M_2)-2=r(H\cap E(M_1))+r(H\cap E(M_2)).$$

Now, for $1\le i\le 2$, $H\cap E(M_i)$ is a proper flat of $M_i$, so that $r(H\cap E(M_i))\le r(M_i)-1$. 
Therefore, $r(H\cap E(M_1))=r(M_1)-1$ and $r(H\cap E(M_2))=r(M_2)-1$.   We conclude that $C_i^*=E(M_i)-H$ is a cocircuit of  $M_i$  and $C^*=C_1^*\cup C_2^*$  such that $C_1^*\cap T=C_2^*\cap T$,  which has exactly two elements.  Note that the converse of the above arguments is also true, thus the proof of the lemma is complete. 
\qed

The following corollary might be of independent interest.

\begin{corollary}\label{3-sum-cocircuit}

Let $M_1$ and $M_2$ be a  binary matroids and $M=M_1\oplus_3M_2$ such that $M_1$ and $M_2$ have the common triangle $T$. Then the following are true:

(i) any cocircuit $C^*$ of $M$  is either a cocircuit of $M_1$ or $M_2$ avoiding $T$,  or  $C^*=C_1^*\Delta C_2^*$ where  $C_i^*$ is a cocircuit of $M_i$ ($i=1,2$) such that $C_1^*\cap T=C_2^*\cap T$, which has exactly two elements.  

(ii) if $C^*$ is either a cocircuit of $M_1$ or $M_2$ avoiding $T$, then $C^*$ is also a cocircuit of $M$. Moreover,  suppose that   $C_i^*$ is a cocircuit of $M_i$ such that $C_1^*\cap T=C_2^*\cap T$, which has exactly two elements.  Then either $C_1^*\Delta C_2^*$ is a cocircuit of $M$, or $C_1^*\Delta C_2^*$ is a disjoint union of two cocircuits $R^*$ and $Q^*$ of  $M$, where $R^*$  and $Q^*$  meet both $M_1$ and $M_2$.

\end{corollary}

\proof  As $M=M_1\oplus_3M_2=P_T(M_1, M_2)\backslash T$, the cocircuits of $M$ are the minimal non-empty members of the set  ${\cal F}=\{D-T$:  $D$ is a cocircuit of $P_T(M_1, M_2)$\}.  If $C^*$ is a cocircuit of $M$,  then $C^*=D-T$ for some cocircuit $D$ of  $P_T(M_1, M_2)$.  By the last lemma, either
(a) $D$ is a cocircuit of $M_1$ or $M_2$ avoiding $T$, or (b) $D=C_1^*\cup C_2^*$ where $C_i^*$ is a cocircuit of $M_i$ ($i=1,2$) such that $C_1^*\cap T=C_2^*\cap T$, which has exactly two elements.  In (a), $C^*=D$, and in (b),  $C^*=C_1^*\Delta C_2^*$. Hence either (i) or (ii) holds in the lemma. 

Conversely, if $C^*$ is either a cocircuit of $M_1$ or $M_2$ avoiding $T$, then clearly $C^*$ is also a cocircuit of $M$, as $C^*=C^*-T$ is clearly a non-empty minimal member of the set 
$\cal F$.   Now suppose that   $C_i^*$ ($i=1,2$) is a cocircuit of $M_i$ such that $C_1^*\cap T=C_2^*\cap T$, which has exactly two elements.  If $C_1^*\Delta C_2^*$ is not a cocircuit of $M$, then it contains a cocircuit $R^*$ of $M$ which is a proper subset of $C_1^*\Delta C_2^*$. Clearly, $R^*$ must meet both $C_1^*$ and $C_2^*$. By (i), $R^*=R_1^*\Delta R_2^*$, 
where $R_i^*$ is a cocircuit of $M_i$ ($i=1,2$) such that $R_1^*\cap T=R_2^*\cap T$, which has exactly two elements. Suppose that $C_1^*\cap T=C_2^*\cap T=\{x, y\}$, then 
$R_1^*\cap T=R_2^*\cap T=\{x, z\}$ or $\{y, z\}$, say the former.  Moreover,  $R_i^*\backslash T$ is a proper subset of $C_i^*\backslash T$ for $i=1,2$ as $T$ does not contain any 
cocircuit of either $M_1$ or $M_2$.  As both $M_1$ and $M_2$ are binary, $Q_i^*=C_i^*\Delta R_i^*$ ($i=1,2$) contains, and indeed, is a cocircuit of $M_i$ such that $Q_1^*\cap T=Q_2^*\cap T=\{y, z\}$. Now it is straightforward to see that $Q_1^*\Delta Q_2^*$ is a minimal non-empty member of $\cal F$ and thus is a cocircuit of $M$. As $C^*=R^*\cup Q^*$, (ii) holds. \qed

The 3-sum of two cographic matroids may not be cographic. However, the following is true. 

\begin{lemma}\label{3sum-cographic}

Suppose that $M_1=M^*(G_1)$ and $M_2=M^*(G_2)$ are both cographic matroids with $u$ and $v$ being vertices of degree three in $G_1$ and $G_2$, respectively.  Label both $uu_i$ and $vv_i$ as $e_i$ ($1\le i\le 3$) so that $T=E(M_1)\cap E(M_2)=\{e_1, e_2, e_3\}$ is the common triangle of $M_1$ and $M_2$.  Then $P_T(M_1,  M_2)=M^*(G)$, where $G$ is obtained by adding a matching $\{u_1v_1, u_2v_2, u_3v_3\}$ between $G_1-u$ and $G_2-v$.  In particular, $M^*(G_1)\oplus_3 M^*(G_2)=M^*(G/e, f , g)$ is also cographic. 
\end{lemma}

\proof We need only show that $P_T(M_1,  M_2)$ and $M^*(G)$ have the same set of cocircuits.  By Lemma \ref{cocircuit}, $C^*$ is a cocircuit of $M=P_T(M_1,  M_2)$ if and only if one of
 the following is true:

(i) $C^*$ is a cocircuit of $M_1$ or $M_2$ avoiding $T$. In other words,  $C^*$ is either a circuit of $G_1$ or a circuit of $G_2$ which does not meet $T$ (i.e., $C^*$ is a circuit of either $G_1-u$ or a circuit of $G_2-v$);

(ii) $C^*=C_1^*\cup C_2^*$ where $C_i^*$ is a cocircuit of $M_i$ such that $C_1^*\cap T=C_2^*\cap T$, which has exactly two elements. In other words, $C^*=C_1^*\cup C_2^*$  where $C_i^*$ ($i=1,2$) is a circuit of $G_i$ containing $u$ and $v$ respectively,  such that $C_1^*\cap T=C_2^*\cap T$, which contains exactly two edges. Now it is easily seen that the set of cocircuits of $M$ is exactly equal to the set of circuits of $M(G)$ (or the set of cocircuits of $M^*(G)$).  In particular, $M^*(G_1)\oplus_3 M^*(G_2)=P_T(M^*(G_1), M^*(G_2))\backslash T=M^*(G)\backslash T=M^*(G/e, f , g)$ is cographic. This completes the proof of the lemma. \qed

The following consequence of the last lemma will be used frequently in the paper. 

\begin{corollary}\label{k3n}

Suppose that $M^*(K_{3, m}), M^*(K^{\prime}_{3, m}), M^*(K_{3, n})\in \cal K$ ($m, n\ge 2$). Then the following are true:

(i) $M^*(K_{3, m}) \oplus_3 M^*(K_{3, n}))\cong M^*(K_{3, m+n-2})$;

(ii) $M^*(K^{\prime}_{3, m}) \oplus_3 M^*(K_{3, n})\cong M^*(K_{3, m+n-2}^{\prime})$;

(iii) $P(M^*(K_{3, m}), M(K_4))$ is cographic and is isomorphic to $M^*(G)$ where $G$ is obtained by putting a 3-edge matching between the 3-partite set of $K_{3, m-1}$ and  the three vertices of $K_{3}$. 

(iv) $M^*(K_{3, m})\oplus_3 M(K_4^\prime))\cong M^*(K_{3, m}^{\prime})$ where $K_4^\prime$ is obtained from $K_4$ by adding a parallel edge to an element in the common triangle $T$ used in the 3-sum. 

(v)  if $M_1\cong M^*(K_{3, m}^{\prime})$, and $M_2\cong M(K_4^\prime)$, then depending on which element in $T$ is in a parallel pair in $M(K_4^\prime)$ and which extra edge was added to $K_{3, m}^{\prime}$ from $K_{3, m}$, the matroid $M_1\oplus_3 M_2$ is either isomorphic to $ M^*(K_{3, m}^{\prime\prime})$ or $M^*(G)$, where $G$ is obtained from $K'_{3,m}$ by adding an edge parallel to the extra edge.  

(vi) if  $M_1\in \cal K$ and $M_2\cong M(K_4^\prime)$, then either $M_1\oplus_3 M_2\in \cal K$  or $M_1\oplus_3 M_2\cong M^*(G)$, where $G$ has a parallel pair which does not meet any triad of $G$. 

(vii) if  $M_1\in \cal K$ and $M_2\in \cal K$, then either $M_1\oplus_3 M_2\in \cal K$  or $M_1\oplus_3 M_2\cong M^*(G)$, where $G$ has at least one parallel pair which does not meet
any triad of $G$. 
\end{corollary}

\proof (i)-(v) are direct consequences of Lemma \ref{3sum-cographic}. Suppose that $M_1\in \cal K$ and is isomorphic to $M^*(K_{3, m})$, $M^*(K'_{3, m})$, $ M^*(K^{''}_{3, m})$, or $ M^*(K^{'''}_{3, m})$. Then either $M_1\oplus_3 M_2\cong M^*(K'_{3, m}), M^*(K^{''}_{3, m})$ or $M^*(K^{'''}_{3, m})$ and thus is in $\cal K$ (in this case, $M_1$ is not isomorphic to $ M^*(K^{'''}_{3, m})$), or  isomorphic to $M^*(G)$, where $G$ is obtained from $K'_{3, m},  K''_{3, m}$, or $K^{'''}_{3, m}$ by adding an edge in parallel to an existing edge added between two vertices of the 3-partite set of $K_{3, m}$. Clearly, this parallel pair does not meet  any triad of $G$. We omit the straightforward and similar proof of (vii).\qed

\begin{corollary}\label{3-sum-starfish}

Let $M$ be a  binary matroid and $M=M_1\oplus_3M_2$ where $M_1$ is a starfish. Suppose that  $M_2$ is a starfish,  or $M_2\cong M(K_4^{\prime})$, or $M_2\cong M^*(G)\in \cal K$: $G\cong K_{3, n}$, $K_{3, n}^{\prime}$, $K_{3, n}^{\prime\prime}$,  or $ K_{3, n}^{\prime\prime\prime}$ ($n\ge 2$). Then  either $M$ is also a starfish, or $M$ has a 2-element cocircuit which does not meet any triangle of $M$. 
\end{corollary}

\proof  Suppose that the starfish $M_1$ uses $s$ Fano matroids and $M_2$ uses $t$ Fano matroids where $s\ge 1$ and $t\ge 0$. Clearly,
in the starfish $M_1$, any triangle is a triad in the corresponding 3-connected graph $G_1\cong K_{3, m}$, $K_{3, m}^{\prime}$, $K_{3, m}^{\prime\prime}$,  or 
$ K_{3, m}^{\prime\prime\prime}$ ($m\ge 2$) used to construct $M_1$. We assume that first $s=1$ and $t=0$. 
Then by the definition of the starfish, $M_1\cong F_7\oplus_3 N_1$, where $N_1\cong M^*(G_1)$,  and either $M_2\cong M(K_4^{\prime})$, or
$M_2\cong M^*(G)$; $G$ is 3-connected  where $G\cong K_{3, n}$, $K_{3, n}^{\prime}$, $K_{3, n}^{\prime\prime}$,  or $ K_{3, n}^{\prime\prime\prime}$ ($n\ge 2$). By Lemma \ref{order}, we have 
that $M=(F_7\oplus_3 N_1)\oplus_3M_2\cong F_7\oplus_3 (N_1 \oplus_3M_2)$ (the condition of the lemma is clearly satisfied). By  Corollary \ref{k3n}, we deduce that either  $N_1 \oplus_3 M_2\in {\cal K}$, or it has a 2-element cocircuit  avoiding any triangle of $N_1 \oplus_3 M_2$. In the former case,  we conclude that $M$ is a starfish. 
In the latter case,  by Corollary \ref{3-sum-cocircuit}, $M$ has a 2-element cocircuit 
avoiding any triangle of $M$.   The general case follows from  an easy induction argument using Lemmas \ref{order} and Corollaries \ref{3-sum-cocircuit} and \ref{k3n}. \qed

\begin{lemma}\label{k4prime}

Suppose that $M\cong M^*(G)$ for a 3-connected graph $G\cong K_{3, n}$, $K_{3, n}^{\prime}$, $K_{3, n}^{\prime\prime}$,  or $ K_{3, n}^{\prime\prime\prime}$ ($n\ge 2$), or $M$ is a starfish. Then for any triangle $T$ of $M$, there are at least two elements $e_1, e_2$ of $T$, such that for each $e_i$ ($i=1, 2$), there is a rooted $K_4^{\prime}$-minor using both $T$ and $e_i$ such that $e_i$ is in a parallel pair. 
\end{lemma}

\proof  Suppose that $M\cong M^*(G)$ for a 3-connected graph $G\cong K_{3, n}$, $K_{3, n}^{\prime}$, $K_{3, n}^{\prime\prime}$,  or $ K_{3, n}^{\prime\prime\prime}$ ($n\ge 2$). When $n\ge 3$, the proof is straightforward. When $n=2$, then $G\cong W_4$ or $K_5\backslash e$, and the result is also true.

Now suppose that $M$ is a starfish constructed by starting from $N\cong  M^*(G)$ for a 3-connected graph $G\cong K_{3, n}$, $K_{3, n}^{\prime}$, $K_{3, n}^{\prime\prime}$,  or $ K_{3, n}^{\prime\prime\prime}$ ($n\ge 2$) with $t$ ($1\le t\le n$) copies of $F_7$ by performing 3-sum operations. Choose an element $f_i$ of $E(M)$ in each copy of 
$F_7$ ($1\le i\le t$).  By the definition of a starfish, and by using Lemma \ref{brylawski}(iii),(v), $M/f_1, f_2, \ldots f_t$ is isomorphic to $N$ containing $T$. Now the result follows from the above paragraph. 
\qed

We will need the following result \cite[11.1]{seymour1980}.

\begin{lemma} \label{graphic}
Let $e$ be an edge of a simple 3-connected graph $G$ on more than four
vertices.
Then either $G\backslash e$ is obtained from a simple 3-connected graph by
subdividing edges or $G/e$ is obtained from a simple 3-connected graph by
adding parallel edges.
\end{lemma}

Let $G=(V,E)$ be a graph and let $x,y$ be distinct elements of $V\cup E$.
By {\it adding an edge between $x,y$} we obtain a graph $G'$ defined as
follows. If $x$ and $y$ are both in $V$, we assume $xy\not\in E$ and we
define $G'=(V,E\cup\{xy\})$; if $x$ is in $V$ and $y=y_1y_2$ is in $E$, we
assume $x\not\in\{y_1,y_2\}$ and we define $G'=(V\cup\{z\},
(E\backslash\{y\})\cup\{xz, y_1z, y_2z\})$; if $x=x_1x_2$ and $y=y_1y_2$
are both in $E$, we define $G'=(V\cup\{u,v\},
(E\backslash\{x,y\})\cup\{ux_1,ux_2,uv, vy_1, vy_2\})$

\begin{lemma} \label{graphic}
Let $G$ be a simple 3-connected graph with a specified triangle $T$. Then
$G$ has a rooted $K_4''$-minor unless $G$ is $K_4$, $W_4$, or Prism.
\end{lemma}

Proof. Suppose the Lemma is false. We choose a counterexample $G=(V,E)$
with $|E|$ as small as possible. Let $x,y,z$ be the vertices of $T$. We
first prove that $G-\{x,y,z\}$ has at least one edge.

Suppose $G-\{x,y,z\}$ is edgeless. Since $G$ is 3-connected, every vertex
in $V-\{x,y,z\}$ must be adjacent to all three of $x,y,z$, which means
that $G=K_{3,n}'''$ for a positive integer $n$. Since $G$ is a
counterexample, $G$ cannot be $K_4$ and thus $G$ contains $K_{3,2}'''$,
which contains a rooted $K_4''$-minor. This contradicts the choice of $G$
and thus $G-\{x,y,z\}$ has at least one edge.

Let $e=uv$ be an edge of $G-\{x,y,z\}$. By Lemma \ref{graphic}, there
exists a simple 3-connected graph $H$ such that at least one of the
following holds:\medskip\\
\indent Case 1. $G\backslash e$ is obtained from $H$ by subdividing edges; \\
\indent Case 2. $G/e$ is obtained from $H$ by adding parallel edges.
\medskip\\
Since $H$ is a proper minor of $G$ and $H$ still contains $T$,  by the minimality of $G$, $H$ has to
be $K_4$, $W_4$, or Prism, because otherwise $H$ (and $G$ as well) would
have a rooted $K_4''$-minor. Now we need to deduce a contradiction in Case
1 and Case 2 for each $H\in\{K_4,W_4$, {\it Prism}\}.

Let $P^+$ be obtained from Prism by adding an edge between two nonadjacent
vertices. Before we start checking the cases we make a simple observation:
with respect to any of its triangles, $P^+$ has a rooted $K_4''$-minor. As
a result, $G$ cannot contain a {\it rooted $P^+$-minor}: a $P^+$-minor in
which $T$ remains a triangle.

We first consider Case 1. Note that $G$ is obtained from $H$ by adding an
edge between some $\alpha,\beta\in V\cup E$. By the choice of $e$, we must
have $\alpha,\beta\not\in V(T)\cup E(T)$. If $H=K_4$ then $G=Prism$, which
contradicts the choice of $G$. If $G=W_4$ or $Prism$, then it is
straightforward to verify that $G$ contains a rooted $P^+$-minor (by
contracting at most two edges), which is a contradiction by the above
observation.

Next, we consider Case 2. Let $w$ be the new vertex created by contracting
$e$. Then $G/e$ is obtained from $H$ by adding parallel edges incident
with $w$. Observe that $w$ has degree three in $H$, for each choice of
$H$. Consequently, as $G$ is simple, $G$ has four, three, or two more
edges than $H$. Suppose $G$ has four or three more edges than $H$. Then
$H$ is $G-u$ or $G-v$. Without loss of generality, let $H=G-u$. Choose
three paths $P_x, P_y, P_z$ in $H$ from $v$ to $x,y,z$, respectively, such
that they are disjoint except for $v$. Now it is not difficult to see that
a rooted $K_4''$-minor of $G$ can be produced from the union of the
triangle $T$, the three paths $P_x,P_y,P_z$, and the star formed by edges
incident with $u$. This contradiction implies that $G$ has exactly two
more edges than $H$. Equivalently, $G$ is obtained from $H$ by adding an
edge between a neighbor $s$ of $w$ and an edge $wt$ with $t\ne s$.

If $H=K_4$ then $G=W_4$, which contradicts the choice of $G$. If $H=W_4$
then $G=W_5$ or $P^+$. In both cases, $G$ contain a rooted $K_4''$-minor,
no matter where the special triangle is. Finally, if $H=Prism$ then $G$
contains a rooted $P^+$-minor, which is impossible by our early
observation. In conclusion, Case 2 does not occur, which completes our
proof. \qed

\begin{lemma} \label{cographic}
Let $M=M^*(G)$ be a 3-connected cographic  matroid with a specified triangle $T$. Then
$M$ has a rooted $K_4''$-minor using $T$ unless $G\cong K_{3, n}$, $K_{3, n}^{\prime}$, $K_{3, n}^{\prime\prime}$,  or $ K_{3, n}^{\prime\prime\prime}$ for some $n\ge 1$.
In particular, if $M^*(G)$ is not graphic, then $n\ge 3$.
\end{lemma}

\proof   Suppose that $M$ does not contain rooted $K_4''$-minor using $T$.  Note that $M^*(G)$ does not have a rooted $K_4''$-minor using $T$ if and only if $G$ does not have a minor obtained from $K_4$ (where $T$ is cocircuit) by subdividing two edges of $T$.  Now we show that $T$ is a vertex triad (which corresponds to a star of degree three). Otherwise, let $G-E(T)=X\cup Y$, where $T$ is a 3-element edge-cut but not a vertex triad. If $G\cong Prism$, then clearly $M^*(G)$ has a rooted $K_4''$-minor; a contradiction.  If $G$ is not isomorphic to a Prism, we can choose a cycle in one side and a vertex in another side which is not incident with any edge of T. Then we can get a rooted $K_4''$-minor; a contradiction again.  Hence the edges of $T$ are all incident to a common vertex $v$ of degree three with neighbors $v_1, v_2$, and $v_3$. A rooted $K_4''$-minor using $T$ exists if and only if $G$ has a cycle missing $v$ and at least two of $v_1, v_2$, and $v_3$.  Hence every cycle of $G-v$ contains at least two of $v_1, v_2$, and $v_3$, and thus  $G-v-v_i-v_j$ is a tree for $1\le i\not=j\le 3$.  Moreover,
$G-v-v_1-v_2-v_3$ has to be an independent set. Otherwise, it is a forest. Take two pedants in a tree, each of which has at least two neighbors in $v_1, v_2$, or $v_3$.  Thus $G-v$ contains a cycle missing at least two vertices of $v_1, v_2,$ and $v_3$. This contradiction shows that $G-v-v_1-v_2-v_3$ is an independent set and thus $G$ is $K_{3, n}$, $K_{3, n}^{\prime}$, $K_{3, n}^{\prime\prime}$,  or $ K_{3, n}^{\prime\prime\prime}$ for some $n\ge 1$. In particular, if $M^*(G)$ is not graphic, then $n\ge 3$. \qed

\begin{lemma}\label{graphic-cographic}

Let $M$ be a 3-connected binary $P_9$-free matroid and $M=M_1\oplus_3 M_2$ where $M_1$ is non-regular, and $M_1$ and $M_2$ have the common triangle $T$. Then

(i) if $M_2$ is graphic, then either $M_2\cong M(G)$ where $G$ is $W_4$ or the Prism, or $M_2\cong M(K_4^{'})$ where $M(K_4^{'})$ is obtained from $M(K_4)$ (which contains $T$) by adding an element parallel to an element of $T$ ; 

(ii) if $M_2$ is cographic but not graphic, then $M_2\cong M^*(G)$, where $G\cong K_{3, n}$, $K_{3, n}^{\prime}$, $K_{3, n}^{\prime\prime}$,  or $ K_{3, n}^{\prime\prime\prime}$ for some $n\ge 3$. 

\end{lemma}

\proof Suppose that $M=P(M_1, M_2)\backslash T$, where $T$ is the common triangle of $M_1$ and $M_2$.  As $M$ is 3-connected,  by \cite[4.3]{seymour1980}, 
both $si(M_1)$ and $si(M_2)$ are 3-connected, and only elements of $T$ can have parallel elements in $M_1$ or $M_2$. Then by Lemma \ref{seymour-F7},  $T$ is contained in a $F_7$-minor in $si(M_1)$.  Now $M_2$ does not contain a rooted $K_4^{''}$-minor using $T$, where $K_4^{''}$ is obtained from this $K_4$ by adding a parallel element to any two of the three elements of $T$ (otherwise, the 3-sum of $M_1$ and $M_2$ contains a $P_9$-minor). 

 If $M_2$ is graphic, then by Lemma \ref{graphic}, $si(M_2)\cong M(G)$ where $G$ is either $W_3, W_4$ or the Prism.  When $G$ is either $ W_4$ or the Prism, then it is easily seen that $M_2$ has to be simple, and thus $M_2\cong W_4$ or Prism. If $G\cong W_3$, then as $M$ is $P_9$-free and $M_2$ has at least seven elements (from the definition of 3-sum),  it is easily seen that $M_2\cong M(K_4^{'})$.  
 
 If $M_2$ is cographic but not graphic, then by Lemma \ref{cographic}, $si(M_2)\cong M^*(G)$, where $G$ is $K_{3, n}$, $K_{3, n}^{\prime}$,
  $K_{3, n}^{\prime\prime}$,  or $ K_{3, n}^{\prime\prime\prime}$ for some $n\ge 3$. If $M_2$ is not simple, then it is straightforward to find a rooted $M(K_4^{''})$-minor using $T$ in $M_2$, thus a $P_9$-minor in $M$; a contradiction. This completes the proof of the lemma. 
 \qed

\begin{lemma}\label{regular1}

Let $M$ be a 3-connected regular matroid with at least six elements and $T$ be a triangle of $M$. Then $M$ has no rooted $M(K_4^{\prime\prime})$-minor using $T$ if and only if 
$M$ is isomorphic to a 3-connected matroid $M^*(K_{3, n})$, $M^*(K_{3, n}^{\prime})$, $M^*(K_{3, n}^{\prime\prime})$,  $M^*(K_{3, n}^{\prime\prime\prime})$ for some $n\ge 1$.
\end{lemma}

\proof 
If $M$ is isomorphic to a 3-connected matroid $M^*(K_{3, n})$, $M^*(K_{3, n}^{\prime})$, $M^*(K_{3, n}^{\prime\prime})$,  $M^*(K_{3, n}^{\prime\prime\prime})$ ($n\ge 1$), then it is straightforward to check  for any triangle $T$,  $M$ has no rooted $M(K_4^{\prime\prime})$-minor using $T$. 

Conversely, suppose that  $M$ is a 3-connected regular matroid with at least six elements and $T$ is a triangle of $M$, such that $M$ has no rooted $M(K_4^{\prime\prime})$-minor using $T$.  If $M$ is internally 4-connected, then by  Theorem \ref{seymour-3-conn}, $M$ is either graphic, cographic, or is 
isomorphic to $R_{10}$. The result follows from 
Lemmas \ref{graphic} and \ref{cographic},  and the fact that $R_{10}$ is triangle-free. So we may assume that $M$ is not internally 4-connected and has a 3-separation $(X, Y)$ 
where $|X|, |Y|\ge 4$. We may assume  that $|X\cap T|\ge 2$.

Suppose that $Y\cap T$ has exactly one element $e$. Then as $T$ is a triangle, $(X\cup e, Y\backslash e)$ is also a 3-separation. If $|Y|=4$, then $Y-e$
 is a triangle or a triad. Moreover, $r(Y)+r^*(Y)-|Y|=2$. As $M$ is 3-connected and binary, $r(Y), r^*(Y)\ge 3$, and thus $r(Y)=r^*(Y)=3$. If $Y-e$ 
is a triangle, then it is not a triad, and thus  $Y$ contains a cocircuit which contains $e$. This is a contradiction as this cocircuit meets $T$ with exactly one element.
Hence $Y-e$ is a triad, and from $r(Y)=3$, there is an element $f\in T,  f\not=e$ such that $Y-f$ is a triangle. In other words, $Y$ forms a 4-element fan. We conclude that 
$M\cong M_1\oplus_3 M(K_4^{\prime})$ by \cite[2.9]{seymour1980} where $S$ is the common triangle of $M_1$ and $M(K_4^{\prime})$, and $M(K_4^{\prime})$ is obtained from $M(K_4)$ (containing $T$) by adding an element $e_1$ in parallel to an element $e$ of $S$.   By switching the label of $e_1$ to $e$ in $M_1$, we obtain a matroid $M'_1$  $(\cong M_1)$ which is isomorphic to a minor of $M$ having triangle $T$.  By \cite[4.3]{seymour1980}, $si(M_1)$ is 3-connected.  Hence by induction, $si(M_1)$ is isomorphic to a 3-connected matroid $M^*(K_{3, m})$, 
$M^*(K_{3, m}^{\prime})$, $M^*(K_{3, m}^{\prime\prime})$,  $M^*(K_{3, m}^{\prime\prime\prime})$  for some $m\ge 1$. As $M$ has no rooted $M(K_4^{\prime\prime})$-minor 
using $T$, we have that $ r_{M_1}(S\cup T)>2$. Moreover,  the element $e_1$ is in two triangles of $si(M_1)$, so $m\le 3$. Now using Lemma \ref{3sum-cographic}, 
it is straightforward to verify that $M\cong W_4\cong M^*(K_{3, 2}^{\prime\prime})$ and thus the Lemma holds. 
Hence we may assume that $|Y|\ge 5$ and thus $|Y\backslash e|\ge 4$.

 Therefore we may assume that $M$ has a separation $(X, Y)$ such that  $T\subseteq X$, 
and both $X$ and $Y$ have at least four elements. 
Hence by \cite[(2.9)]{seymour1980},  $M=M_1 \oplus_3 M_2$ where $M_1$ and $M_2$ are isomorphic to minors of $M$ having the common triangle $S$, and $T$ is a triangle of $M_1$. Moreover,  $|E(M_i)| < |E(M)|$ for $i=1, 2$, and both $si(M_1)$ and $si(M_2)$ are $3$-connected \cite[(4.3)]{seymour1980}.  First assume that each element of $S$ is parallel to an element of $T$ in $M_1$. Then by Lemma \ref{k4}, $si(M_1)$ contains a rooted $M(K_4)$-minor using $T$. As each element of $T$ in $M_1$ is in a parallel pair, we conclude that 
$M$ has a rooted  $M(K_4^{\prime\prime})$-minor using $T$; a contradiction.

So we may assume that at least one element of $T$ is  not parallel to an element of $S$ (as $M$ is binary, there are at least two such elements). As $si(M_1)$ is a 3-connected minor of $M$, it has  no rooted $M(K_4^{\prime\prime})$-minor using $T$. By induction, $si(M_1)\cong M^*(K_{3, s})$, $M^*(K_{3, s}^{\prime})$, $M^*(K_{3, s}^{\prime\prime})$, $M^*(K_{3, s}^{\prime\prime\prime})$ for some $s\ge 2$,  or  $si(M_1)\cong M(K_4)$.   Remove all elements of $M_1$ not in the set $S\cup T$ in $P_S(M_1, M_2)$. Then every element of $T\backslash S$ is parallel to an element of $S\backslash T$. Contracting all elements of $S\backslash T$, we obtained a minor of $M$ isomorphic to $M_2$  and $T$ is a triangle of this minor.  By induction again,  $si(M_2)\cong M^*(K_{3, t})$, $M^*(K_{3, t}^{\prime})$, $M^*(K_{3, t}^{\prime\prime})$, $M^*(K_{3, t}^{\prime\prime\prime})$ for some $t\ge 2$,  or  $si(M_2)\cong M(K_4)$.  Suppose that $si(M_i)\cong M(K_4)$ for some $i=1, 2$. Then  as $M_i$ have at least seven elements and $M$ has no rooted $M(K_4^{\prime\prime})$-minor
 using $T$, we deduce that $M_i\cong M(K_4^{\prime})$. As $M$ has no $M(K_4^{\prime\prime})$-minor containing $T$, and $M$ is 3-connected, using Corollary \ref{k3n},  
 it is routine to verify that  $M \cong M^*(K_{3, n}), M^*(K_{3, n}^{\prime})$, $M^*(K_{3, n}^{\prime\prime})$, or $M^*(K_{3, n}^{\prime\prime\prime})$  for some $n\ge 2$.
\qed

\begin{corollary}\label{regular}

Let $M$ be a $3$-connected binary non-regular $P_9$-free matroid. Suppose that $M=M_1\oplus_3 M_2$  such that $M_1$ and $M_2$ have the 
common triangle $T$. 
 If $M_2$ is regular, then $M_2$ is isomorphic to a 3-connected matroid   $M^*(K_{3, n})$, $M^*(K_{3, n}^{\prime})$, $M^*(K_{3, n}^{\prime\prime})$, or $M^*(K_{3, n}^{\prime\prime\prime})$ ($n\ge 2$), or $M_2\cong M(K_4^{\prime})$ where $M(K_4^{\prime})$ is obtained from $M(K_4)$ (containing $T$) by adding an element in parallel to an element of $T$. 
\end{corollary}

\proof  As $M$ is 3-connected,  by \cite[4.3]{seymour1980}, 
both $si(M_1)$ and $si(M_2)$ are 3-connected, and only elements of $T$ can have parallel elements in $M_1$ or $M_2$. As $M$ is non-regular and $M_2$ is regular, 
$si(M_1)$ is non-regular and thus (by Lemma \ref{seymour-F7}) has a $F_7$-minor containing the common triangle  $T$ of $M_1$ and $M_2$. 
As $M$ is $P_9$-free, $M_2$ has no rooted $M(K_4^{\prime\prime})$-minor using $T$. By Lemma \ref{regular1}, $si(M_2)$ is isomorphic to a 3-connected matroid $M^*(K_{3, n})$, $M^*(K_{3, n}^{\prime})$, $M^*(K_{3, n}^{\prime\prime})$,  $M^*(K_{3, n}^{\prime\prime\prime})$ ($n\ge 2$), or $M(K_4)$. Now using Lemma \ref{k4prime},  it is 
straightforward to check that either $M_2\cong M(K_4^{\prime})$, or $M_2$ is simple, and $M_2\cong  M^*(K_{3, n})$, $M^*(K_{3, n}^{\prime})$, $M^*(K_{3, n}^{\prime\prime})$, or $M^*(K_{3, n}^{\prime\prime\prime})$ ($n\ge 2$).  \qed

Now we are ready to prove our main theorem.

\noi {\it Proof of Theorem \ref{main}}.  Suppose that a starfish $M$ is constructed from a 3-connected cographic 
matroid $N$ by consecutively applying the 3-sum operations with $t$ copies of $F_7$, where $N\cong M^*(G)$;  $G\cong K_{3, n}, K_{3, n}^{\prime}, 
K_{3, n}^{\prime\prime}$, or $K_{3, n}^{\prime\prime\prime}$ for some $n\ge 2$.  First we show that $M$ is 3-connected. We use induction on $t$. When $t=0$, $N$ is 3-connected. Suppose that $M$ is 
3-connected for $t<k\le n$. Now suppose that $t=k$. Then $M=M_1\oplus_3 F$, where $F\cong F_7$ and $M_1$ and $F$ share the common triangle $T$. 
Take an element $f$ of $E(F)\cap E(M)$. Then by Lemma \ref{brylawski}, $M/f=P(M_1, F/e)\backslash T\cong M_1$, which is a starfish with $t=k-1$, and thus is 3-connected by induction. If $M$ is not 3-connected, then $f$ is either in a loop of $M$, or is in a
cocircut of size one or two. Clearly, $M$ does not have any loop, thus $f$ is in a cocircuit $C^*$ of $M$ with size one or two. As $P(M_1, F)$ is 3-connected, it does not contain 
any cocircuit of size less than three.  Hence  $C^*\cup T$ contains a cocircuit $D^*$ of $P(M_1, F)$. As $P(M_1, F)$ is binary, $D^*\cap T$ has exactly two elements, and 
thus $D^*$ has at most four elements. As $T$ contains no cocircuit of either
 $M_1$ or $F$, by Lemma \ref{cocircuit}, $F\cong F_7$ has a cocircuit of size at most three meeting two elements of $T$. This contradiction shows that $M$ is 3-connected.

Next we show that if $M$ is one of the matroid listed in (i)-(iv), then $M$ is $P_9$-free.  By Theorem \ref{maini4c} and the fact that all spikes and regular matroids are $P_9$-free, 
we need only show that any  starfish
is $P_9$-free. We use induction on the number of elements of the starfish $M$.
By the definition, the unique smallest starfish has nine elements, and is isomorphic to $P_9^*$. Clearly, $P_9^*$ is  $P_9$-free. Suppose that any starfish with less than $n$ $(\ge 10)$ elements is  
$P_9$-free. Now suppose that we have a starfish $M$ with $n$ elements. Suppose,  on the contrary, that $M$ has a $P_9$-minor. Then by the Splitter Theorem (Theorem \ref{SST}), 
there is an element $e$ in $M$ such that either $M\backslash e$ or $M/e$ is 3-connected having a $P_9$-minor.  Note that the elements of a starfish consists of 
two types: those are subsets of $E(N)$ (denote this set by $K$), or those are in part of copies of $F_7$ (denote this set by $F$).  Then $E(M)=K\cup F$. First we assume 
that $e\in F$. Then $M=M_1\oplus_3 M_2$, where $M_1$ is either one of $M^*(K_{3, n})$, $M^*(K_{3, n}^{\prime})$, $M^*(K_{3, n}^{\prime\prime})$, or 
$M^*(K_{3, n}^{\prime\prime\prime})$, or a starfish with fewer elements; $M_2\cong F_7$, and $e\in E(M_2)$. By the construction of the starfish and 
Lemma \ref{brylawski},  $M/e\cong M_1$ and is either cographic or a smaller starfish and therefore does not contain a $P_9$-minor; a contradiction.
Therefore $M\backslash e$ is 3-connected and 
contains a $P_9$-minor. But then by Lemma \ref{order}, $M\backslash e\cong P(M_1, M(K_4))\backslash T$. 
By Corollary \ref{k3n}, as $M\backslash e$ is 3-connected, we conclude that $M\backslash e$ is a smaller starfish and therefore is $P_9$-free. This contradiction shows that $e\in K$.

 If $e$ is in a triangle of $M$, then $M/e$ is not 3-connected, and thus $M\backslash e$ is 3-connected and contains a $P_9$-minor. Each triangle of 
 $M$ is corresponding to a triad in $G$. By Lemmas \ref{brylawski} and \ref{order} again, we can do the deletion $N\backslash e$ first, then perform the 3-sum operations 
 with copies of $F_7$.  Note that $N\backslash e\cong M^*(G/e)$ where $G\cong K_{3, n}$, $K_{3, n}^{\prime}$, $K_{3, n}^{\prime\prime}$,  or $ K_{3, n}^{\prime\prime\prime}$
  ($n\ge 2$).   As $M\backslash e$ is 3-connected and thus simple, we deduce that $n\ge 3$, 
$N\cong M^*(K_{3, n})$ or $M^*(K_{3, n}^{\prime})$, and $N\backslash e\cong M^*(K_{3, n-1}^{\prime\prime})$, or $M^*(K_{3, n-1}^{\prime\prime\prime})$.  Therefore, 
$M\backslash e$ is another starfish and does not contain any $P_9$-minor by induction; a contradiction. Finally assume that $e\in K$ is not in any triangle of $M$. Then $e$ is not in any
triad of $G$. Hence  if $n=2$, then $G\cong K'''_{3,2}$. 
 As $G/e$ has parallel elements, the matroid $N\backslash e$ has serial-pairs, 
and thus $M\backslash e$ is not 3-connected, we conclude that $M/e$ is 3-connected having a $P_9$-minor. Note that $N\cong  M^*(K_{3, n}^{\prime})$,  $M^*(K_{3, n}^{\prime\prime})$, or $M^*(K_{3, n}^{\prime\prime\prime})$ ($n\ge 2$), and thus $N/e \cong  M^*(K_{3, n})$, $M^*(K_{3, n}^{\prime})$, or $M^*(K_{3, n}^{\prime\prime})$, which is still 3-connected. We conclude again, by Lemma \ref{brylawski},  that $M/e$ is a smaller starfish than $M$, thus cannot contain any $P_9$-minor. This contradiction completes the proof of the first part.

Now suppose that $M$ is a 3-connected binary matroid with no $P_9$-minor. We may assume that $M$ is not regular. 
 If $M$ is  internally 4-connected, then the theorem follows from 
Theorem  \ref{maini4c}. 
Now suppose that $M$  is neither regular nor internally 4-connected. We show that $M$ is either a spike or a starfish. Suppose that $|E(M)|\le 9$. As $M$ is not internally 4-connected, $M$ is not $F_7$ or $F_7^*$. Hence $|E(M)|\ge 8$. Then $M$ is $AG(3,2),  S_8, Z_4,  Z_4^*$ (all spikes), or $P_9^*$,  which is the 3-sum of $F_7$ and $W_4=M^*(K_{3, 2}^{''})$, thus is a starfish. We conclude that the result holds for $|E(M)|\le 9$. Now suppose that $|E(M)|\ge 10$. As $M$ is not internally 4-connected,  $M=M_1 \oplus_3 M_2=P(M_1, M_2)\backslash T$,  where $M_1$ and $M_2$ are isomorphic to minors of $M$ (\cite[4.1]{seymour1980}) and $T=\{x, y, z\}$ is the common triangle of $M_1$ and $M_2$. Moreover,  $|E(M_i)| < |E(M)|$ for $i=1, 2$, and both $si(M_1)$ and $si(M_2)$ are $3$-connected \cite[(4.3)]{seymour1980}.  The only possible parallel element(s) of either $M_1$ or $M_2$ are those in the common triangle. As $M$ has no $P_9$-minor,  and $M_1$ and $M_2$ are isomorphic to minors of $M$, we deduce that neither $si(M_1)$ nor $si(M_2)$ has a $P_9$-minor.  By induction, the theorem holds for both $si(M_1)$ and $si(M_2)$.  As $M$ is not regular, at least one of $si(M_1)$ and $si(M_2)$, say $si(M_1)$,  is not regular.  

\noindent {\bf Claim}: $M_1$ (and $M_2$) is simple unless both $si(M_1)$ and $si(M_2)$ are spikes. 

Suppose not and we may assume that $x$ in $T$ has a parallel element $x_1$ in $M_1$.   By Lemma \ref{seymour-F7}, $T$ is in a $F_7$-minor of $M_1$ plus a parallel element $x_1$. By induction, $si(M_2)$ is either regular and 3-connected, or one of the $16$  internally 4-connected non-regular minors of $Y_{16}$ (thus is $F_7$ since it has a triangle); or  is a spike or a   starfish. Moreover, $si(M_1)$ is
 either one of the $16$  internally 4-connected non-regular minors of $Y_{16}$ (thus is $F_7$); or  is a spike or a   starfish.  Suppose that $si(M_2)$ is not a spike. Then either $si(M_2)$ is regular or is a  starfish. 
By Lemmas  \ref{k4prime} and \ref{regular},   either $M_2\cong M(K_4^{\prime})$ where $M(K_4^{'})$ is obtained from $M(K_4)$ (which contains $T$) by adding an element parallel to an element of $T$, or $T$ is in a rooted $M(K_4^{\prime})$-minor  of $M_2$ using $T$ (obtained from $M(K_4)$ containing $T$ by adding an element parallel to either $y$ or $z$). In either case, as $M$ is simple, we conclude that $M$ contains a $P_9$-minor, a contradiction. Hence $si(M_2)$ is a spike thus contains an $F_7$-minor containing $T$.  Now if $si(M_1)$ is not a spike, then $si(M_1)$ is a starfish. Again using Lemma  \ref{k4prime}, it is easily checked that $M$ has a $P_9$-minor; a contradiction. 
Therefore $M_1$ is simple unless both $si(M_1)$ and $si(M_2)$ are spikes. A similar argument shows that $M_2$ is also simple unless both $si(M_1)$ and $si(M_2)$ are spikes.

\noindent {\bf Case 1}:  $si(M_2)$ is regular. By Lemma \ref{regular}, $M_2$ is either graphic or cographic. Moreover, 

(i) if $M_2$ is graphic, then either $M_2\cong M(G)$ where $G$ is $W_4$ or the Prism,  or $M_2\cong M(K_4^{\prime})$ 
where $M(K_4^{'})$ is obtained from $M(K_4)$ (which contains $T$) by adding an element parallel to an element of $T$;  and

(ii) if $M_2$ is cographic but not graphic, then $M\cong M^*(G)$, where $G\cong K_{3, n}$, $K_{3, n}^{\prime}$, $K_{3, n}^{\prime\prime}$,  or $ K_{3, n}^{\prime\prime\prime}$ for some 
$n\ge 3$.  

By the above claim, both $M_1$ and $M_2$ are simple.  Moreover, $M_1$ is 3-connected, non-regular, and $P_9$-free. By induction, $M_1$ is either one of the $16$  internally 4-connected non-regular minors of $Y_{16}$ (therefore is $F_7$ as $M_1$ has a triangle); or $M_1$ is a spike or a  starfish. That is,  either $M_1$ is a spike or a  starfish. 
If $M_1$ is a starfish,  by Lemma \ref{3-sum-starfish}, $M=M_1\oplus_3 M_2$ is also a starfish. Thus we may assume that $M_1$ is a spike which contains a triangle.  
Then $M_1$ is either $F_7, S_8$,  
$Z_s$ ($s\ge 4$) or $Z_s\backslash y_s$ for some $s\ge 5$.  Suppose that $M_1$ is $F_7$.  Then $M=F_7\oplus_3 M_2$ is either $S_8$ (not possible as $M$ has at least 10 elements) or a starfish by the definition of a starfish.  Suppose that $M_1$ is $Z_s$ ($s\ge 4$) or $Z_s\backslash y_s$ for some $s\ge 5$ and suppose that $M_2$ is not isomorphic to $M(K_4^{\prime})$. Then $M_1$ has a $Z_4$-restriction containing $T$. Clearly, such restriction contains a $F_7^{\prime}$-minor which is obtained from $F_7$ (which contains $T$) by adding an element parallel to the tip of the spike, say $x$ in $T$. By Lemma  \ref{k4prime}, $T$ is in a $M(K_4^{\prime})$- minor  of $M_2$ which is obtained from $K_4$ containing $T$ by adding an element parallel to an element $z\not=x$ of $T$. Thus we can find  a $P_9$-minor in $M$, a contradiction. Suppose that $M_1$ is $Z_s$ ($s\ge 4$) or $Z_s\backslash y_s$ for some $s\ge 5$ and suppose that $M_2\cong M(K_4^{\prime})$. If the extra element $e$ of $M(K_4^{\prime})$ added to $M(K_4)$ is not parallel to $x$ in $M_2$, then
using the previously mentioned $F'_7$-minor of $M_1$ containing $T$ and the $M(K_4^{\prime})$-minor containing $e$, we obtain a $P_9$-minor of $M$; a contradiction. 
Now it is straightforward to see that  $M\cong Z_{s+2}\backslash y_{s+2}$ ($s\ge 4$) which is a spike, or   $Z_{s+2}\backslash y_s,  y_{s+2}$ ($s\ge 5$). The latter
case does not happen as $\{y_s,  y_{s+2}\}$ would be  a 2-element cocircuit, but $M$ is 3-connected.  Finally we assume that $M_1\cong S_8=F_7\oplus_3 M(K_4^{\prime})$ with tip $x$. Then $M=(F_7\oplus_3 M(K_4^{\prime}))\oplus_3 M_2$.  By Lemma 
 \ref{order}, $M=F_7\oplus_3 (M(K_4^{\prime})\oplus_3 M_2)$. By Corollary \ref{regular}, $M_2$ is isomorphic to a 3-connected cographic matroid $M^*(K_{3, n})$, $M^*(K_{3, n}^{\prime})$, $M^*(K_{3, n}^{\prime\prime})$, or $M^*(K_{3, n}^{\prime\prime\prime})$ ($n\ge 2$), or $M_2\cong M(K_4^{\prime})$.  If  $M_2\cong M(K_4^{\prime})$, then $|E(M)|=9$; a contradiction. Thus $M_2$ is not isomorphic to $M(K_4^{\prime})$.  By Corollary \ref{k3n}, 
$M(K_4^{\prime})\oplus_3 M_2\cong M^*(G)$, where $G\cong K_{3, n}^{\prime}$, $K_{3, n}^{\prime\prime}$,  or $ K_{3, n}^{\prime\prime\prime}$ for some $n\ge 2$, or $M(K_4^{\prime})\oplus_3 M_2$ contains a 2-element cocircuit which does not meet any triangle of $M(K_4^{\prime})\oplus_3 M_2$. In this case, by Corollary \ref{3-sum-cocircuit}, this 2-element cocircut would also be a cocircuit of $M$. As $M$ is 3-connected, we conclude that the latter does not happen, and that $M$ is still a starfish.

\noindent {\bf Case 2}:   Neither $M_1$ nor $M_2$ is regular. By induction and the fact that both $M_1$ and $M_2$ have a triangle, that $si(M_1)$ is either a spike containing a triangle or a starfish, and so is $si(M_2)$.

Case 2.1: Both $si(M_1)$ and $si(M_2)$ are  starfishes. By the above claim, both $M_1$ and $M_2$ must be simple matroids.  Now by Lemma \ref{3-sum-starfish}, $M$ is also a starfish.

Case 2.2:  One of $si(M_1)$  and $si(M_2)$, say the former,  is a spike. Suppose that $si(M_2)$ is a starfish. By the claim, both $M_1$ and $M_2$ are simple.  
As $M_1$ contains the triangle $T$, it is either $Z_s$ ($s\ge 3$) or $Z_s\backslash y_s$ for some $s\ge 4$. If $M_1\cong Z_3\cong F_7$, by the definition of a starfish, $M$ is also a starfish.  If $M_1\cong Z_s$ ($s\ge 4$) or $Z_s\backslash y_s$ for some $s\ge 5$, then $M_1$ contains a $Z_4$ as a restriction which contains $T$. But $Z_4$ contains a $F_7^{\prime}$-minor 
containing $T$ where $F_7^{\prime}$ is obtained from $F_7$ by adding  an element in parallel to the tip $x$ of $M_1$. By Lemma \ref{k4prime}, $T$ is in a $M(K_4^{\prime})$-minor  of 
$M_2$ which is obtained from $M(K_4)$ containing $T$ by adding an element parallel to $y$ or $z$. We conclude that $M$ contains a $P_9$-minor, a contradiction. Now suppose that $M_1\cong Z_4\backslash y_4\cong S_8=F_7\oplus_3 M(K_4^{\prime})$ with tip $x$. Then $M=(F_7\oplus_3 M(K_4^{\prime}))\oplus_3 M_2$. By Lemma
\ref{order}, $M=F_7\oplus_3 (M(K_4^{\prime})\oplus_3 M_2)$. By  Corollary \ref{3-sum-starfish}, $M(K_4^{\prime})\oplus_3 M_2$ is either a starfish, or  $M(K_4^{\prime})\oplus_3 M_2$ and  thus $M$ contains a 2-element cocircuit. As $M$ is 3-connected, we conclude that the latter does not happen, and that $M$ is still a starfish by 
the definition of a starfish.

 Hence we may assume that $si(M_2$) is also a spike.  As $si(M_2)$ contains a triangle also,   it is either  $Z_t$ ($t\ge 3$) or $Z_t\backslash y_t^{\prime}$ for some $t\ge 4$. 
 Suppose that $si(M_1)$ and $si(M_2)$ do not share a common tip, say $si(M_1)$ has tip $x$ and $si(M_2)$ has tip $z$. Then neither matroid is isomorphic to $F_7$ as any 
 element of $T$ can be considered as a tip then. We first assume  either $si(M_1)$ or $si(M_2)$, say $si(M_1)$,  has at least nine elements. Then $M_1$ has a $Z_4$-restriction 
 containing $T$, thus has a $F_7^{\prime}$-minor (with a parallel pair containing $x$) containing $T$.  The matroid  $si(M_2)$ has a $S_8$-restriction, thus has a 
 $M(K_4^{\prime})$-minor (with a parallel pair containing $z$) containing $T$. By Lemma \ref{brylawski}, we conclude that $M$ has a $P_9$-minor;  a contradiction. 
 Hence both $si(M_1)$ and $si(M_2)$ have exactly eight elements and both are isomorphic to $S_8$. Now if either $M_1$ or $M_2$ is not simple, then similar to the argument 
 above, one can get a $P_9$-minor; a contradiction. Hence both matroid are simple. Now it is straightforward to see that $M\cong F_7\oplus_3 W_4\oplus_3 F_7$, which is a starfish. 
 
 Therefore we may assume that $si(M_1)$ and $si(M_2)$ share a common tip, say $x$. First assume that a non-tip element in $T$, say $y$,  is in a parallel pair of either $M_1$ or $M_2$, say $M_1$.  As $M$ is both simple and $P_9$-free, it is easily seen that $M_2$ has to be simple. Since any element of $T$ can be considered as a tip in $F_7$, we deduce that both $si(M_1)$ and $M_2$ have at least 8 elements.  If one of these two matroids has at least 9 elements, then it contains a $Z_4$-restriction containing $T$. Such a restriction contains a $F_7^{\prime}$-minor containing $T$ with $x$ being in a parallel pair. At the same time, $si(M_i)$ contains a $M(K_4)$-minor 
 containing $T$ for $i=1, 2$. Noting that $y$ is in a parallel pair of $M_1$, we deduce that $M$ contains a $P_9$-minor; a contradiction. Hence we may assume that 
 both  $si(M_1)$ and $M_2$ contain exactly 8 elements. Now it is easily seen that $M_1$ contains a $F_7^{\prime}$-minor containing $T$ with $y$ being in a parallel pair. At the same time, $si(M_2)$ contains a $M(K_4^{\prime})$-minor containing $T$ with $x$ being in a parallel pair. This is a contradiction as $M$ now contains a $P_9$-minor.

 So from now on we may assume that if $M_1$ or $M_2$ is not simple, then only $x$ could be in a parallel pair. Indeed, as $M$ is simple, at most one of $M_1$ and
  $M_2$ is not simple. Suppose that one of $M_1$ and $M_2$, say $M_1$,  is not simple, then either $M\cong Z_{s+t}$, $M\cong Z_{s+t}\backslash y_s$, 
  $M\cong Z_{s+t}\backslash y_t^{\prime}$, or $M\cong Z_{s+t}\backslash y_s, y_t^{\prime}$, all of which are spikes except the last matroid. The last matroid, $M\cong Z_{s+t}\backslash y_s, y_t^{\prime}$, however, contains a cocircuit $\{y_s, y_t^{\prime}\}$, contradicting to the fact that $M$ is 3-connected. 
Finally assume that both $M_1$ and $M_2$ are simple. Then $M\cong Z_{s+t}\backslash x$, $M\cong Z_{s+t}\backslash x, y_s$, $M\cong Z_{s+t}\backslash x, y_t^{\prime}$, or $M\cong Z_{s+t}\backslash x, y_s, y_t^{\prime}$, all of which are spikes except the last matroid.   The last matroid, $M\cong Z_{s+t}\backslash x, y_s, y_t^{\prime}$, again, contains a cocircuit $\{y_s, y_t^{\prime}\}$; a contradiction. This completes the proof of Case 2.2, thus the proof of the theorem. \qed

\section*{Acknowledgements}
G. Ding's research is partially supported by  NSA grant
H98230-14-1-0108.

\end{document}